\documentclass[reqno]{amsart}
\usepackage{amsmath, amssymb, amsthm, epsfig}
\usepackage{hyperref, latexsym}
\usepackage{url} 
\def\today{\ifcase\month\or
  January\or February\or March\or April\or May\or June\or
  July\or August\or September\or October\or November\or December\fi
  \space\number\day, \number\year}

 \newtheorem{theorem}{Theorem}
 \newtheorem{lemma}[theorem]{Lemma}
 
 \newtheorem{corollary}[theorem]{Corollary}
 \theoremstyle{definition}

 \theoremstyle{remark}

 \newcommand{\C}{\mathbb{C}}
 \newcommand{\R}{\mathbb{R}}

 \newcommand{\Z}{\mathbb{Z}}
 
 \newcommand{\tF}{\widehat{F}}

 \newcommand{\uh}{\widehat{h}}

 \newcommand{\tf}{\widehat{f}}

 \newcommand{\dt}{\text{\rm d}t}

 \newcommand{\dr}{\text{\rm d}r}
 
 \newcommand{\dx}{\text{\rm d}x}
 \newcommand{\dy}{\text{\rm d}y}

 \newcommand{\domega}{\text{\rm d}\omega}
 \newcommand{\dtau}{\text{\rm d}\tau}
 \newcommand{\dsigma}{\text{\rm d}\sigma}

 \newcommand{\ot}{\otimes}
 \newcommand{\la}{\langle}
 \newcommand{\ra}{\rangle}
 \newcommand{\deta}{\text{\rm d}\eta}
 \newcommand{\dzeta}{\text{\rm d}\zeta}
 \newcommand{\dxi}{\text{\rm d}\xi}

\begin{document}

\title[Strichartz inequality]{A sharp inequality for the Strichartz norm}
\author[E. Carneiro]{Emanuel Carneiro}
\thanks{The author acknowledges support from CAPES/FULBRIGHT grant BEX 1710-04-4.}
\date{\today}
\address{Department of Mathematics, University of Texas at Austin, Austin, TX 78712-1082.}
\email{ecarneiro@math.utexas.edu}
\allowdisplaybreaks
\numberwithin{equation}{section}
\begin{abstract} 
Let $u:\R \times \R^n \to \C$ be the solution of the linear Schr\"odinger equation
\begin{equation*}
\left\{ 
\begin{array}{ccl}
iu_t  + \Delta u & =&  0 \\
u(0,x) &=& f(x).
\end{array}
\right.
\end{equation*}
In the first part of this paper we obtain a sharp inequality for the Strichartz norm $\|u(t,x)\|_{L^{2k}_tL^{2k}_x(\R \times\R^n)}$, where $k\in \Z$, $k \geq 2$ and $(n,k) \neq (1,2)$, that admits only Gaussian maximizers. As corollaries we obtain sharp forms of the classical Strichartz inequalities in low dimensions (works of Foschi \cite{F} and Hundertmark - Zharnitsky \cite{HZ}) and also sharp forms of some Sobolev-Strichartz inequalities. In the second part of the paper we express Foschi's \cite{F} sharp inequalities for the Schr\"odinger and wave equations in the broader setting of sharp restriction/extension estimates for the paraboloid and the cone. 
\end{abstract}
\maketitle
\section{Introduction}
Let $u:\R \times \R^n \to \C$ be the solution of the linear  Schr\"odinger equation
\begin{equation}\label{Intro1}
\left\{ 
\begin{array}{ccl}
iu_t + \Delta u & =& 0 \\
u(0,x) &=& f(x)\,.
\end{array}
\right.
\end{equation}
The homogeneous Strichartz estimates (see \cite{C}) are inequalities of the type
\begin{equation}\label{Intro2}
\|u(t,x)\|_{L^q_tL^r_x(\R \times\R^n)} \leq C \|f\|_{L^2(\R^n)},
\end{equation}
with
\begin{equation*}
\|u(t,x)\|_{L^q_tL^r_x(\R \times\R^n)} = \left[\int_{\R} \left( \int_{\R^n}|u(t,x)|^r \dx\right)^{q/r} \dt \,\right]^{1/q}.
\end{equation*}
The pair of exponents $(q,r)$ is admissible if 
\begin{equation*}
\frac{2}{q} + \frac{n}{r} = \frac{n}{2}\,,
\end{equation*}
with $2\leq q,r\leq \infty$ and  $(q,r,n)\neq (2,\infty,2)$. The sharp forms of the Strichartz inequalities were first investigated in a paper by Kunze \cite{K}, who showed the existence of maximizers in the case $n=1$, $(q,r)=(6,6)$, by concentration-compactness techniques. Later, Foschi \cite{F} and Hundertmark-Zharnitsky \cite{HZ} independently obtained the sharp constants in the cases $n=1$, $(q,r) = (6,6)$; and $n=2$, $(q,r) = (4,4)$; showing that the only maximizers are Gaussians. They conjectured that in the case $q=r= 2 + 4/n$, $n\geq 3$, the extremals for the Strichartz inequalities should be given by Gaussians. Recently, Shao \cite{Sh} showed that maximizers do exist for the non-endpoint Strichartz inequalities ($q\neq 2$ if $n\geq 3$ and $q\neq 4$ if $n=1$) in all dimensions.

In this note we generalize the beautiful argument of \cite{HZ} to prove the following sharp inequality for the Strichartz norm.

\begin{theorem}\label{thm1}
Let $u:\R \times \R^n \to \C$ be the solution of the linear Schr\"odinger equation {\rm (\ref{Intro1})}. For $k \in \Z$, $k \geq 2$ and $(n,k) \neq (1,2)$ we have
\begin{equation}\label{main ineq}
\|u(t,x)\|_{L^{2k}_tL^{2k}_x(\R \times\R^n)} \leq \left( C_{n,k} \int_{\R^{nk}} |\tF(\eta)|^2 \,K(\eta)^{\tfrac{n(k-1) - 2}{2}}\, \deta \right)^{1/2k},
\end{equation}
with 
\begin{equation}\label{constant}
C_{n,k} = \left[2^{n(k-1)-1}\,k^{n/2}\, \pi^{(n(k-1)-2)/2}\, \Gamma\left(\tfrac{n(k-1)}{2}\right)\right]^{-1}.
\end{equation}
On the right hand side of {\rm (\ref{main ineq})} we write $\eta \in \R^{nk}$ as $\eta = (\eta_1,\eta_2,...,\eta_k)$ with each $\eta_i \in \R^n$; $F(\eta) = f(\eta_1)f(\eta_2)...f(\eta_k)$; and the kernel 
\begin{equation*}
K(\eta) = \frac{1}{k}\sum_{1\leq i < j \leq k} |\eta_i - \eta_j|^2.
\end{equation*}
This inequality is sharp and equality occurs if and only if $f$ is a Gaussian.
\end{theorem}

Throughout this paper we will adopt the definition of the Fourier transform of the function $f:\R^n \to \C$ given by
\begin{equation*}
\tf(\omega) = \frac{1}{(2\pi)^{n/2}} \int_{\R^n} e^{-i\omega \cdot x} f(x) \,\dx.
\end{equation*}
We observe that the solution of (\ref{Intro1}) can be given in terms of the Fourier transform 
\begin{equation}\label{solution}
u(t,x) = \frac{1}{(2\pi)^{n/2}} \int_{\R^n} e^{i x \cdot \omega}\, e^{-it|\omega|^2}\, \tf(\omega)\, \domega.
\end{equation}
The maximizers in Theorem \ref{thm1} should be understood in the following way: if $\tf$ is a measurable function such that the right hand side of {\rm (\ref{main ineq})} is finite, and equality occurs in (\ref{main ineq}), then $\tf$ must be a Gaussian, and so is $f$. Here we shall always refer as Gaussians the functions of the form
\begin{equation}\label{special gaussian}
f(x) = e^{A|x|^2 + b\cdot x + C}\,,
\end{equation}
where $A, C \in \C$, $b \in \C^n$ and $\Re(A)<0$. The term $A$ is the covariance of the Gaussian $f$.

Some interesting inequalities arise from Theorem \ref{thm1}. First, we present the sharp forms of the classical Strichartz inequalities in low dimensions.
\begin{corollary} In dimension $n=1$ we have
\begin{equation}\label{Intro3}
\|u(t,x)\|_{L^6_tL^6_x(\R \times\R)} \leq 12^{-1/12} \|f\|_{L^2(\R)},
\end{equation}
and 
\begin{equation}\label{Intro4}
\|u(t,x)\|_{L^8_tL^4_x(\R \times\R)} \leq 2^{-1/4}\|f\|_{L^2(\R)}.
\end{equation}
In dimension $n=2$ we have
\begin{equation}\label{Intro5}
\|u(t,x)\|_{L^4_tL^4_x(\R \times\R^2)} \leq 2^{-1/2} \|f\|_{L^2(\R^2)}.
\end{equation}
These inequalities are sharp and equality occurs if and only if $f$ is a Gaussian.
\end{corollary}
The sharp forms (\ref{Intro3}) and (\ref{Intro5}) are the ones discovered by Foschi \cite{F} and Hundertmark-Zharnitsky \cite{HZ}. They are a direct consequence of Theorem \ref{thm1}. The novelty here is (\ref{Intro4}), which is obtained
by taking $f(x,y) = g(x)g(y)$ in (\ref{Intro5}) and exploiting the product structure of the problem. It is interesting to notice the persistence of the Gaussian maximizers in a case where $q \neq r$.

By using the fact that 
\begin{equation}\label{gen}
\int_{\R^n \times \R^n} g(x)\,g(y)\,x\cdot y \,\dx \,\dy \geq 0\,,
\end{equation}
for any real valued function $g$, with equality for example if $g$ is radial, one obtains some sharp Sobolev-Strichartz inequalities in low dimensions.
\begin{corollary} In dimension $n=1$ we have
\begin{equation}\label{Intro6}
 \|u(t,x)\|_{L^{10}_tL^{10}_x(\R \times\R)} \leq  (2\sqrt{5}\pi)^{-1/10} \|f'\|_{L^2(\R)}^{1/5} \|f\|_{L^2(\R)}^{4/5},
\end{equation}
\begin{equation}\label{Intro7}
 \|u(t,x)\|_{L^{12}_tL^6_x(\R \times\R)} \leq  (6\pi)^{-1/12} \|f'\|_{L^2(\R)}^{1/6} \|f\|_{L^2(\R)}^{5/6},
\end{equation}
and 
\begin{equation}\label{Intro8.3}
 \|u(t,x)\|_{L^{16}_tL^4_x(\R \times\R)} \leq  (8\pi)^{-1/16} \|f'\|_{L^2(\R)}^{1/8} \|f\|_{L^2(\R)}^{7/8}.
\end{equation}
In dimension $n=2$ we have
\begin{equation}\label{Intro8}
 \|u(t,x)\|_{L^6_tL^6_x(\R \times\R^2)} \leq  (12\pi)^{-1/6} \|\nabla f\|_{L^2(\R^2)}^{1/3} \|f\|_{L^2(\R^2)}^{2/3},
\end{equation}
and 
\begin{equation}\label{Intro8.2}
 \|u(t,x)\|_{L^8_tL^4_x(\R \times\R^2)} \leq  (16\pi)^{-1/8} \|\nabla f\|_{L^2(\R^2)}^{1/4} \|f\|_{L^2(\R^2)}^{3/4}.
\end{equation}
In dimension $n=4$ we have
 \begin{equation}\label{Intro8.1}
 \|u(t,x)\|_{L^4_tL^4_x(\R \times\R^4)} \leq  (32\pi)^{-1/4} \|\nabla f\|_{L^2(\R^4)}^{1/2} \|f\|_{L^2(\R^4)}^{1/2}.
\end{equation}
These inequalities are sharp and equality occurs if and only if $f$ is a Gaussian.
\end{corollary}
Inequalities (\ref{Intro6}), (\ref{Intro8}) and (\ref{Intro8.1}) follow directly from Theorem \ref{thm1} and (\ref{gen}). To obtain (\ref{Intro7}) and (\ref{Intro8.3}) one should put $f(x,y) = g(x)g(y)$ in (\ref{Intro8}) and (\ref{Intro8.2}), respectively, and exploit the product structure. In an analogous manner one obtains (\ref{Intro8.2}) by putting $f(x,y,z,k) = g(x,y)g(z,k)$ in (\ref{Intro8.1}).

\subsection{Sharp restriction/extension estimates}
It has been known for a long time the equivalence of decay inequalities for the space-time norm of the solutions of certain evolution equations and restriction estimates for the Fourier transform over curved surfaces. The classical reference on the subject is Strichartz original paper \cite{S}, but seminal ideas can already be observed in the work of H\"{o}rmander \cite[Corollary 1.3]{H}. 

The Schr\"odinger and wave equations are related to the restriction problem for the paraboloid and cone, respectively,
\begin{equation}\label{parab}
S_{parab} := \{ (\tau, \omega) \in \R \times \R^{n}\,:\, \tau = |\omega|^2\}\,,
\end{equation}
and 
\begin{equation}\label{cone}
S_{cone} := \{ (\tau, \omega) \in \R \times \R^{n}\,:\, \tau = |\omega|\}\,.
\end{equation}
We endow these surfaces $S \subset \R^{n+1}$ with canonical measures $\dsigma$ given by
\begin{equation}\label {meas parab}
\int_{S_{parab}}g(\tau, \omega) \,\dsigma = \int_{\R^n} g(|\omega|^2, \omega)\, \domega\,,
\end{equation}
and
\begin{equation}\label{meas cone}
\int_{S_{cone}}g(\tau, \omega) \,\dsigma = \int_{\R^n} g(|\omega|, \omega)\, \frac{\domega}{|\omega|}\,.
\end{equation}
In this setting, the restriction estimates are a priori inequalities of the form 
\begin{equation}\label{rest}
\|\uh|_S\|_{L^{p'}(S;\,\dsigma)} \leq C_{p,q,S}\|h\|_{L^{q'}(\R^{n+1})}.
\end{equation}
The scaling invariance tells us that the global estimate (\ref{rest}) can only hold for $p' = nq/(n+2)$ in the case of the paraboloid and $p' = (n-1)q/(n+1)$ in the case of the cone. On the other hand, Knapp's example shows that we must have $q >(2n+2)/n$ for the paraboloid and $q > 2n/(n-1)$ for the cone. The restriction conjecture asserts that these are sufficient conditions in each case for (\ref{rest}) to hold, and so far it has been proved for the range $q > (2n+6)/(n+1)$ in both cases, the paraboloid by Tao \cite{T2} and the cone by Wolff \cite{W}. We refer the reader to \cite{T} for a survey on the recent progress on the restriction conjecture.

A duality argument using Parseval's identity shows that 
\begin{align}\label{duality}
\begin{split}
C_{p,q,S} & = \sup_{\|h\|_{L^{q'}(\R^{n+1})} =1} \ \ \|\uh|_S\|_{L^{p'}(S;\,\dsigma)} \\
& = \sup_{\|h\|_{L^{q'}(\R^{n+1})} =1}\ \ \sup_{\|g\|_{L^{p}(S;\,\dsigma)} = 1} \ \ \left|\int_S \uh(\tau, \omega)\, g(\tau,\omega)\,\dsigma\right| \\
& = \sup_{\|g\|_{L^{p}(S;\,\dsigma)} = 1}\ \ \sup_{\|h\|_{L^{q'}(\R^{n+1})} =1} \ \ \left|\int_{\R^{n+1}} h(t,x)\, \widehat{g\dsigma}\,(t,x)\,\dt \,\dx\right| \\
& = \sup_{\|g\|_{L^{p}(S;\,\dsigma)} = 1}\ \ \|\widehat{g\dsigma}\|_{L^q(\R^{n+1})}\,.
\end{split}
\end{align}
Therefore (\ref{rest}) is equivalent to the extension estimate
\begin{equation}\label{ext}
\|\widehat{g\dsigma}\|_{L^{q}(\R^{n+1})} \leq C_{p,q,S} \|g\|_{L^{p}(S;\,\dsigma)}\,,
\end{equation}
for all smooth functions $g$ on $S$, where $\widehat{g\dsigma}$ is the Fourier transform of the measure $g\dsigma$:
\begin{align*}
\widehat{(g\dsigma)}(t,x) :=\frac{1}{(2\pi)^{(n+1)/2}} &\int_{S} g(\tau, \omega) \,e^{-i(t\tau + \omega\cdot x)} \,\dsigma.
\end{align*}
In the case of the paraboloid, from (\ref{solution}) we see that the solution of the Schr\"odinger equation (\ref{Intro1}) satisfies
\begin{equation*}
u(t,-x) = (2\pi)^{1/2}\widehat{g\dsigma}(t,x)\,,
\end{equation*}
with $g(|\omega|^2, \omega) = \tf(\omega)$. Therefore, (\ref{ext}) is equivalent to the inequality
\begin{equation}\label{Strich}
\|u(t,x)\|_{L^{q}_tL^{q}_x(\R \times\R^n)} \leq (2\pi)^{1/2}C_{p,q,S} \|\tf\|_{L^{p}(\R^n)}.
\end{equation}
From the equivalence of (\ref{rest}), (\ref{ext}) and (\ref{Strich}), the sharp forms (\ref{Intro3}) and (\ref{Intro5}) discovered by Foschi \cite{F} and Hundertmark-Zharnitsky \cite{HZ} immediately translate into sharp restriction/extension estimates for the paraboloid.
\begin{theorem}
Let $S$ be the paraboloid defined in {\rm(\ref{parab})} endowed with the measure $\dsigma$ defined in {\rm(\ref{meas parab})}. We have
\begin{equation}\label{parab1}
\|\widehat{g\dsigma}\|_{L^{6}(\R^{2})} \leq (2\pi)^{-1/2} 12^{-1/12} \|g\|_{L^{2}(S;\,\dsigma)}\,,
\end{equation}
and 
\begin{equation}\label{parab2}
\|\widehat{g\dsigma}\|_{L^{4}(\R^{3})} \leq (4\pi)^{-1/2} \|g\|_{L^{2}(S;\,\dsigma)}.
\end{equation}
These inequalities are sharp. Equality occurs in (\ref{parab1}) and (\ref{parab2}) if and only if
\begin{equation}\label{max parab}
g(|\omega|^2, \omega) = e^{A|\omega|^2 + b\cdot \omega + C}\,,
\end{equation}
where $A, C \in \C$, $b \in \C^n$ and $\Re(A)<0$.
\end{theorem}
For simplicity, we presented above the sharp extension inequality. One can deduce the dual sharp restriction inequality (\ref{rest}) for the paraboloid and find the maximizing functions $h(t,x)$ by using the condition for equality in the duality argument (\ref{duality}) (H\"{o}lder's inequality) 
\begin{equation}\label{duality condition}
h = C\,|\widehat{g\dsigma}|^{\tfrac{q}{q'} -1}\, \overline{\widehat{g\dsigma}}\,,
\end{equation}
for a complex constant $C$ and $g$ given by (\ref{max parab}).

In the same spirit, sharp restriction/extension inequalities for the cone are implicit in Foschi's work \cite{F} for the wave equation.

\begin{theorem}\label{thm5}
Let $S$ be the cone defined in {\rm(\ref{cone})} endowed with the measure $\dsigma$ defined in {\rm(\ref{meas cone})}. We have
\begin{equation}\label{cone1}
\|\widehat{g\dsigma}\|_{L^{6}(\R^{3})} \leq (2\pi)^{1/3} \|g\|_{L^{2}(S;\,\dsigma)}\,,
\end{equation}
and 
\begin{equation}\label{cone2}
\|\widehat{g\dsigma}\|_{L^{4}(\R^{4})} \leq (2\pi)^{1/4} \|g\|_{L^{2}(S;\,\dsigma)}.
\end{equation}
These inequalities are sharp. Equality occurs in (\ref{cone1}) and (\ref{cone2}) if and only if
\begin{equation}\label{max cone}
g(|\omega|, \omega) = e^{A|\omega| + b\cdot \omega + C}\,,
\end{equation}
where $A, C \in \C$, $b \in \C^n$ and $|\Re(b)| < -\Re(A)$.
\end{theorem}

We will give a brief proof of Theorem \ref{thm5} in section 4, indicating the basic changes that have to be made in Foschi's argument. Again, the maximizers $h(t,x)$ for the dual restriction inequalities (\ref{rest}) can be obtained from the duality condition (\ref{duality condition}) with $g$ given by (\ref{max cone}). It would be a very interesting line of research to investigate other sharp constants in the broader setting of restriction/extension estimates and to understand the role that the special functions (\ref{max parab}) and (\ref{max cone}) play in these inequalities.\\

We shall see in this paper that the natural generalization of the argument of Hundertmark-Zharnitsky \cite{HZ} leads to the inequality in Theorem \ref{thm1}, which maintains the Gaussian maximizers, but is weaker than (\ref{Strich}). Indeed, one can show that for 
\begin{equation*}
q = 2k  \ \ \ \textrm{and} \ \ \ p = \frac{2nk}{2nk - n - 2}\,,
\end{equation*}
the following inequality holds
\begin{equation}\label{weak}
\|\tf\|_{L^{p}(\R^n)} \leq C  \left(\int_{\R^{nk}} |\tF(\eta)|^2 \,K(\eta)^{\tfrac{n(k-1) - 2}{2}}\, \deta \right)^{1/2k}\,.
\end{equation}
This is a consequence of the following three inequalities:
\begin{enumerate}
\item[(i)] A basic inequality for real numbers:
\begin{equation*}
K(\eta)^{\tfrac{n(k-1) - 2}{2}} \geq C \sum_{1\leq i < j \leq n} |\eta_i - \eta_j|^{n(k-1) - 2}\,;
\end{equation*}
\item[(ii)] The reversed Hardy-Littlewood-Sobolev inequality due to W. Beckner \cite{B}: 
\begin{equation*}
\int_{\R^n \times \R^n} |g(x)||x-y|^{\lambda}|h(y)|\,\dx\,\dy \geq C(n,\lambda)\, \|g\|_{L^{\tfrac{2n}{2n+\lambda}}(\R^n)}\|h\|_{L^{\tfrac{2n}{2n+\lambda}}(\R^n)}\,,
\end{equation*}
where $\lambda > 0$, the sharp constant given by
\begin{equation*}
C(n,\lambda) = \pi^{\lambda/2} \frac{\Gamma(n/2 + \lambda/2)}{\Gamma(n + \lambda/2)} \left[\frac{\Gamma(n)}{\Gamma(n/2)}\right]^{1 + \lambda/n},
\end{equation*}
and the only maximizers being $g(x) = c\,h(x)$, $c\in \C$ a constant, and 
\begin{equation*}
h(x) = A (B^2 + |x-x_0|^2)^{-(2n+\lambda)/2}\,,
\end{equation*}
for some $A \in \C$, $0 \neq B \in \R$ and $x_0 \in \R^n$. For our purposes it suffices to use this inequality in the following format
\begin{equation*}
\int_{\R^n \times \R^n} |\tf(\eta_i)|^2|\tf(\eta_j)|^2 |\eta_i - \eta_j|^{n(k-1) - 2}\, \deta_i\, \deta_j \geq C \|\tf\|_{L^r(\R^n)}^4\,,
\end{equation*}
where $r = 4n/(n(k+1)-2)$;\\
\item[(iii)] H\"{o}lder's inequality: 
\begin{equation*}
\|\tf\|_{L^{p}(\R^n)}^{2k} \leq \|\tf\|_{L^r(\R^n)}^4\|\tf\|_{L^2(\R^n)}^{2k-4}.
\end{equation*}
\end{enumerate}
Inequality (\ref{weak}) will be used later in section 3.
\section{Proof of Theorem \ref{thm1} - the sharp inequality}

The proof of Theorem \ref{thm1} given here follows closely the outline of Hundertmark and Zharnitsky \cite{HZ}. As we are interested in an a priori estimate, in this section we suppose that $f \in C^{\infty}_0(\R^n)$. Throughout the proof of Theorem \ref{thm1} we reserve the variables $\eta$ and $\xi$ to be in $\R^{nk}$ and write $\eta = (\eta_1,\eta_2,...,\eta_k)$ with each $\eta_i \in \R^n$. We have also defined $F(\eta) = f(\eta_1)f(\eta_2)...f(\eta_k)$ and $K(\eta) = \frac{1}{k}\sum_{1\leq i < j \leq k} |\eta_i - \eta_j|^2$. Let us write 
\begin{equation*}
F_1(\eta) = \tF(\eta)K(\eta)^{\tfrac{n(k-1) -2}{4}}.
\end{equation*}
In the space $L^2(\R^{nk})$, let $E$ be the closed subspace consisting of the functions invariant under any orthonormal transformation (rotation here for short) $R$ that fixes the vectors $\alpha_1, \alpha_2,..., \alpha_n \in \R^{nk}$ given by
\begin{equation}\label{S2.1}
\alpha_i = (e_i, e_i,...,e_i)  \ \ (k  \ \ \textrm{times}),
\end{equation}
where $e_i = (0,0,...,1,...,0)$ is the $i$-th canonical vector in $\R^n$. Denote by $P_E: L^2(\R^{nk}) \to L^2(\R^{nk})$ the orthogonal projection operator onto the subspace $E$. The heart of the matter is the following representation lemma.

\begin{lemma}[Representation Lemma]\label{lem4}
Let $u:\R \times \R^n \to \C$ be the solution of the Schr\"odinger equation \rm{(\ref{Intro1})}. Then
\begin{equation*}
\int_{\R \times \R^n} |u(t,x)|^{2k}\, \dx \,\dt = C_{n,k} \,\langle P_E(F_1), F_1\ra_{L^2(\R^{nk})}.
\end{equation*}
with the constant $C_{n,k}$ defined in (\ref{constant}).
\end{lemma}
\begin{proof}
Using the representation (\ref{solution}) for the solution $u(t,x)$ we obtain
\begin{equation*}
|u(t,x)|^{2k} = \frac{1}{(2\pi)^{nk}}\int_{\R^{nk} \times \R^{nk}} e^{i x \cdot (\sum \eta_i  - \sum \xi_i)}\, e^{-it(|\eta|^2 - |\xi|^2)}\, \tF(\eta)\overline{\tF(\xi)}\, \deta\, \dxi,
\end{equation*}
where $\eta = (\eta_1,\eta_2,...,\eta_k)$ and $\xi = (\xi_1,\xi_2,...,\xi_k)$, with each $\eta_i$ and $\xi_i$ in  $\R^n$. Integrating with respect to $x$ and $t$ and using that, as distributions, the n-dimensional delta function $\delta_n(w) = (2\pi)^{-n} \int_{\R^n} e^{-ix\cdot w} \dx$\,, one arrives at
\begin{align*}
\int_{\R \times \R^n}& |u(t,x)|^{2k}\, \dx \,\dt \\
& = \frac{1}{(2\pi)^{n(k-1) -1}}\int_{\R^{nk} \times \R^{nk}} \delta_n\left(\sum_{i=1}^k \eta_i  - \sum_{i=1}^k \xi_i\right)\, \delta\bigl(|\eta|^2 - |\xi|^2\bigr)\, \tF(\eta)\overline{\tF(\xi)}\, \deta\, \dxi\\
& = \frac{1}{(2\pi)^{n(k-1) -1}}\int_{\R^{nk} \times \R^{nk}} \left(\prod_{i=1}^{n} \delta\bigl((\eta - \xi)\cdot \alpha_i\bigr)\right)\, \delta\bigl(|\eta|^2 - |\xi|^2\bigr)\, \tF(\eta)\overline{\tF(\xi)}\, \deta\, \dxi.
\end{align*}
We will rewrite the last equation in the following strategic way
\begin{align*}
&\int_{\R \times \R^n} |u(t,x)|^{2k}\, \dx \,\dt \\
& = \frac{1}{(2\pi)^{n(k-1) -1}}\int_{\R^{nk} \times \R^{nk}} \frac{\left(\prod_{i=1}^{n} \delta\bigl((\eta - \xi)\cdot \alpha_i\bigr)\right)\, \delta\bigl(|\eta|^2 - |\xi|^2\bigr)}{\bigl(K(\eta)K(\xi)\bigr)^{\tfrac{n(k-1) -2}{4}}}\, F_1(\eta)\overline{F_1(\xi)}\, \deta\, \dxi.
\end{align*}
The insight now is to recognize the last expression as a quadratic form associated to a self-adjoint operator. Indeed, for $G \in C^{\infty}_0(\R^{nk})$ define the operator
\begin{equation}\label{def of A}
AG(\xi) = \frac{1}{(2\pi)^{n(k-1) -1}}\int_{\R^{nk}} \frac{\left(\prod_{i=1}^{n} \delta\bigl((\eta - \xi)\cdot \alpha_i\bigr)\right)\, \delta\bigl(|\eta|^2 - |\xi|^2\bigr)}{\bigl(K(\eta)K(\xi)\bigr)^{\tfrac{n(k-1) -2}{4}}}\, G(\eta)\, \deta\,.
\end{equation}
In this context we have 
\begin{equation*}
\int_{\R \times \R^n} |u(t,x)|^{2k}\, \dx \,\dt = \la AF_1,F_1\ra_{L^2(\R^{nk})}\,.
\end{equation*}
Our objective is to show that the operator $A$ is a multiple of the projection operator $P_E$. We start by showing that $A$ is a bounded operator in $L^2(\R^{nk})$, via the following lemma.

\begin{lemma}\label{lem5}
\begin{enumerate}
\item[(i)] For all $\xi \in \R^{nk}$ the measure
\begin{equation*}
m_{\xi}(\deta) = \frac{k^{n/2} \, \Gamma\left(\tfrac{n(k-1)}{2}\right)}{\pi^{n(k-1)/2}} \frac{\left(\prod_{i=1}^{n} \delta\bigl((\eta - \xi)\cdot \alpha_i\bigr)\right)\, \delta\bigl(|\eta|^2 - |\xi|^2\bigr)}{\bigl(K(\eta)K(\xi)\bigr)^{\tfrac{n(k-1) -2}{4}}}\, \deta
\end{equation*}
is a probability measure on $\R^{nk}$.\\
\item[(ii)] For all Borel measurable sets $B \subset \R^{nk}$, we have
\begin{equation*}
\int_{\R^{nk}} m_{\xi}(B) \,\dxi = |B|\,,
\end{equation*} 
where $|B|$ denotes the Lebesgue measure of $B$.
\end{enumerate}

\end{lemma}
\begin{proof}Throughout this proof let us write 
\begin{equation*}
C = \frac{k^{n/2} \,\Gamma\left(\tfrac{n(k-1)}{2}\right)}{\pi^{n(k-1)/2}}\,.
\end{equation*}
Observe that in the support of the delta functions we have $\sum \eta_i = \sum \xi_i$ and $|\eta|^2 = |\xi|^2$. This implies that $K(\eta) = K(\xi)$, since
\begin{align}\label{K}
\begin{split}
K(\eta) &= \frac{1}{k}\sum_{1\leq i < j \leq k} |\eta_i - \eta_j|^2 \\
& = |\eta|^2 - \frac{|\eta_1 + \eta_2 + ...+ \eta_k|^2}{k} =  |\eta|^2 - \frac{\sum_{i=1}^n (\eta \cdot \alpha_i)^2}{k}\,.
\end{split}
\end{align}
Therefore we have
\begin{equation}\label{prob measure}
m_{\xi}(\R^{nk}) = \frac{C}{K(\xi)^{\tfrac{n(k-1)-2}{2}}} \int_{\R^{nk}}\left(\prod_{i=1}^{n} \delta\bigl((\eta - \xi)\cdot \alpha_i\bigr)\right)\, \delta\bigl(|\eta|^2 - |\xi|^2\bigr)\, \deta\,.
\end{equation}
Let $\{\tilde{e}_j\}$, $1 \leq j \leq nk$, be the canonical vectors in $\R^{nk}$. Change the variable $\eta$ in the integration (\ref{prob measure}) by a rotation $R$ that sends $\alpha_i$ to $\sqrt{k}\tilde{e}_i$ for $1\leq i \leq n$. We obtain
\begin{align*}
m_{\xi}&(\R^{nk}) = \frac{C}{K(\xi)^{\tfrac{n(k-1)-2}{2}}} \int_{\R^{nk}} \delta_n\left(\sqrt{k}\eta_1 - \sum\xi_i\right)\, \delta\bigl(|\eta|^2 - |\xi|^2\bigr)\, \deta\\
& = \frac{C}{k^{n/2}K(\xi)^{\tfrac{n(k-1)-2}{2}}}\int_{\R^{n(k-1)}} \delta\left(\sum_{i=2}^{k}|\eta_i|^2 - K(\xi)\right)\, \deta_2 \,\deta_3...\deta_k\\
& = \frac{C\,\bigl|S^{n(k-1)-1}\bigr|}{k^{n/2}K(\xi)^{\tfrac{n(k-1)-2}{2}}}\int_0^{\infty}\delta(r^2 - K(\xi))\, r^{n(k-1)-1} \dr\\
\\
& = \frac{C\,\bigl|S^{n(k-1)-1}\bigr|}{2\,k^{n/2}K(\xi)^{\tfrac{n(k-1)-2}{2}}}\int_0^{\infty}\delta(t - K(\xi))\, t^{\tfrac{n(k-1)-2}{2}} \dt \\
\\
& = \frac{C\,\bigl|S^{n(k-1)-1}\bigr|}{2\,k^{n/2}} = 1\,,
\end{align*} 
and this proves (i). To prove (ii), just observe the symmetry of the measure $m$ with respect to the variables $\eta$ and $\xi$,
\begin{align*}
\int_{\R^{nk}} m_{\xi}(B) \,\dxi &= \int_{\R^{nk}} \int_B  \frac{C\,\left(\prod_{i=1}^{n} \delta\bigl((\eta - \xi)\cdot \alpha_i\bigr)\right)\, \delta\bigl(|\eta|^2 - |\xi|^2\bigr)}{\bigl(K(\eta)K(\xi)\bigr)^{\tfrac{n(k-1) -2}{4}}}\, \deta\, \dxi\\
& =\int_B \int_{\R^{nk}}\frac{C\,\left(\prod_{i=1}^{n} \delta\bigl((\eta - \xi)\cdot \alpha_i\bigr)\right)\, \delta\bigl(|\eta|^2 - |\xi|^2\bigr)}{\bigl(K(\eta)K(\xi)\bigr)^{\tfrac{n(k-1) -2}{4}}}\, \dxi\, \deta\\
& = \int_B m_{\eta}(\R^{nk})\, \deta = \int_B \deta = |B|.
\end{align*}
\end{proof}
We now return to the proof of the Representation Lemma \ref{lem4}. Note the the operator $A$ can be written as
\begin{equation*}
AG(\xi) = C_{n,k}\int_{\R^{nk}} G(\eta) m_{\xi}(\deta)\,.
\end{equation*}
The boundedness of the operator $A$ in $L^2(\R^{nk})$ follows from an application of Lemma \ref{lem5} and Jensen's inequality
\begin{align*}
\|AG&\|_{L^2(\R^{nk})}^2 = C_{n,k}^{2} \int_{\R^{nk}}\left|\int_{\R^{nk}} G(\eta) m_{\xi}(\deta)\right|^2 \dxi \leq C_{n,k}^{2} \int_{\R^{nk}}\int_{\R^{nk}} |G(\eta)|^2 m_{\xi}(\deta)\, \dxi \\
& = C_{n,k}^{2}  \int_{\R^{nk}} |G(\eta)|^2 \int_{\R^{nk}} m_{\xi}(\deta)\,\dxi = C_{n,k}^{2} \int_{\R^{nk}} |G(\eta)|^2 \deta = C_{n,k}^{2} \|G\|_{L^2(\R^{nk})}^2.
\end{align*}
We thus arrive at 
\begin{equation*}
 \|AG\|_{L^2(\R^{nk})} \leq C_{n,k}\|G\|_{L^2(\R^{nk})}, 
\end{equation*}
proving that the operator $A$ extends to a bounded operator from $L^2(\R^{nk})$ to $L^2(\R^{nk})$. It remains to show that $A$ is a multiple of the projection operator $P_E$. Let $R$ be a rotation on $\R^{nk}$ fixing the vectors $\alpha_1, ...,\alpha_n$. It is clear from (\ref{def of A}) and (\ref{K}) that 
\begin{equation*}
 AG(R\xi) = AG(\xi),
\end{equation*}
therefore $A$ maps $L^2(\R^{nk})$ into the subspace $E$. From the fact that the operator $A$ is self-adjoint we can show that $A(E^{\perp}) = 0$. It remains to prove that $A$ acts like a multiple of the identity on $E$. For this, consider a function $H \in C^{\infty}_0(\R \times \R \times  ...\times \R \times \R^+)$ and write
\begin{equation}\label{form of E}
G(\eta) = H(\eta\cdot \alpha_1, \eta \cdot \alpha_2, ..., \eta \cdot \alpha_n, |\eta|^2).
\end{equation}
Certainly $G$ is a function in $E$, and from definition (\ref{def of A}) we find that, for a $G$ of the form (\ref{form of E}),
\begin{equation*}
AG(\xi) = C_{n,k}G(\xi).
\end{equation*}
Since the functions of the form (\ref{form of E}) are dense in $E$, we conclude that $A = C_{n,k}I$ on $E$. We have proved that $A = C_{n,k} P_E$ and this concludes the lemma.
\end{proof}

The proof of the inequality proposed in Theorem \ref{thm1} is then a trivial consequence of the Representation Lemma \ref{lem4}. In fact,
\begin{align}\label{ineq thm1}
\begin{split}
 \int_{\R \times \R^n} |u(t,x)|^{2k}\, \dx \,\dt& = C_{n,k} \langle P_E(F_1), F_1\ra_{L^2(\R^{nk})} \leq C_{n,k} \|F_1\|_{L^2(\R^{nk})}^2\\
& = C_{n,k} \int_{\R^{nk}} |\tF(\eta)|^2 \,K(\eta)^{\tfrac{n(k-1) - 2}{2}}\, \deta\,.
\end{split}
\end{align}
It remains to investigate when equality in (\ref{ineq thm1}) can be attained. A necessary and sufficient condition is that the function $F_1(x)$ belongs to the subspace $E$.

\section{Proof of Theorem \ref{thm1} - Gaussian maximizers}

We investigate here under which conditions the function
\begin{equation*}
F_1(\eta) = \tF(\eta)K(\eta)^{\tfrac{n(k-1) -2}{4}}
\end{equation*}
belongs to the subspace $E$. Let us say that a measurable function $G:\R^{nk} \to \C$ satisfies the property ($\star$) if $G$ is invariant under all the rotations $R$ that fix the vectors $\alpha_1, \alpha_2,..., \alpha_n$. In this setting, $G \in E$ if and only if $G \in L^2(\R^{nk})$ and satisfies ($\star$).

From (\ref{K}) we see that $K(x)$ satisfies ($\star$). Therefore, we must have $\tF(\eta) = \tf(\eta_1)\tf(\eta_2)...\tf(\eta_k)$ satisfying ($\star$), and we shall prove that under these symmetries $\tf$ must be a Gaussian. The proof will be divided in five steps.\\
\\
{\bf Step 1.}  Let $g:\R^n \to \C$ be a measurable function such that $G(\eta) = g(\eta_1)g(\eta_2)...g(\eta_k)$ satisfies
\begin{equation}\label{cond}
\int_{\R^{nk}} |G(\eta)|^2 \,K(\eta)^{\tfrac{n(k-1) - 2}{2}}\, \deta < \infty\,.
\end{equation}
Then $g \in L^p(\R^n)$ for $p = \tfrac{2nk}{2nk - n - 2}$.\\

This was proved in (\ref{weak}). From now on we fix $p = \tfrac{2nk}{2nk - n - 2}$.\\
\\
{\bf Step 2.} Let $g \in L^p(\R^n)$ be such that $G(\eta)$ satisfies the property ($\star$). Then $g$ is a product of one-dimensional functions.\\

We shall write here each $\eta_i \in \R^n$ as $\eta_i = (\eta_{i1}, \eta_{i2},...,\eta_{in})$. If $g \in L^p(\R^n)$ is nonzero, there exists a cube $J = \prod_{i=1}^n [a_i,b_i] \subset \R^n$ such that 
\begin{equation*}
 \int_J g(y)\, \dy = A \neq 0.
\end{equation*}
Consider the orthonormal transformation $R$ in $\R^{nk}$ that simply switches the coordinates $\eta_{11}$ and $\eta_{21}$ on $\eta = (\eta_1,...,\eta_k)$. Naturally, this transformation fixes the vectors $\alpha_i$ and thus the relation  $G(Rx) = G(x)$ implies
\begin{align}\label{eq for g}
\begin{split}
g(\eta_{11}, \eta_{12},...,\eta_{1n})&g(\eta_{21}, \eta_{22},...,\eta_{2n})g(\eta_3)...g(\eta_k) \\
& = g(\eta_{21}, \eta_{12},...,\eta_{1n})g(\eta_{11}, \eta_{22},...,\eta_{2n})g(\eta_3)...g(\eta_k).
\end{split}
\end{align}
Integrating both sides of (\ref{eq for g}) with respect to $\deta_2\deta_3...\deta_k$ on $J \times J \times...\times J$ we find that 
\begin{align}\label{eq for g 2}
\begin{split}
A^{k-1} g(&\eta_{11}, \eta_{12},...,\eta_{1n}) \\
&= A^{k-2} \int_{a_1}^{b_1} g(\eta_{21}, \eta_{12},...,\eta_{1n})\, \deta_{21} \int_{J'}g(\eta_{11}, \eta_{22},...,\eta_{2n})\, \deta_2'\,,
\end{split}
\end{align}
where  $J' = \prod_{i=2}^n [a_i,b_i]$ and $\deta_2' = \deta_{22}\deta_{23}...\deta_{2n}$. Expression (\ref{eq for g 2}) plainly says that 
\begin{equation}\label{eq for g 3}
 g(\eta_{11}, \eta_{12},...,\eta_{1n}) = w_1(\eta_{11})\,h_1(\eta_{12},...,\eta_{1n}).
\end{equation}
By repeating this argument we arrive at 
\begin{equation}\label{eq for g 4}
 g(\eta_{11}, \eta_{12},...,\eta_{1n}) = w_j(\eta_{1j})\,h_j(\eta_{11},...,\eta_{1(j-1)},\eta_{1(j+1)},...,\eta_{1n}),
\end{equation}
for $j =2,...,n$. Expressions (\ref{eq for g 3}) and (\ref{eq for g 4}) are sufficient to conclude that 
\begin{equation*}
g(\eta_{11}, \eta_{12},...,\eta_{1n}) = g_1(\eta_{11})g_2(\eta_{12})...g_n(\eta_{1n}).
\end{equation*}
\\
{\bf Step 3.}  Suppose that all $g_i$'s are smooth and non-vanishing. Then all $g_i$'s are Gaussians with the same covariance. Therefore $g$ is itself a Gaussian.\\

Let $R_{12}$ be a rotation on $\R^{2n}$ fixing the vectors $\beta_i = \frac{1}{\sqrt{2}}(e_i,e_i)$, $i=1,2,...,n$. Observe that the rotation on $\R^{nk}$ given by
\begin{equation}\label{matrix}
R = \left[
\begin{array}{ccccc}
R_{12} & &&&\\
& I & &0&\\
&&I &&\\
&0&&\ddots&\\
&&&& I
\end{array}
\right]
\end{equation}
fixes the vectors $\alpha_i = (e_i,e_i,...,e_i) \in \R^{nk}$. Among all the possible rotations $R$ given by this form, we will choose a simple rotation $R_{12}$ to work with. Let us denote the tensor product $a \ot b$ of two vectors $a= (a_1, a_2,...,a_n)$ and $b = (b_1,b_2,...,b_n)$ in $\R^n$ as the $n \times n$ matrix $[a_ib_j]$, corresponding to the linear transformation $x \mapsto \langle x, b \rangle a$. Consider the orthonormal basis of $\R^{2n}$ formed by the vectors $\beta_i = \frac{1}{\sqrt{2}}(e_i,e_i)$ and $\gamma_i = \frac{1}{\sqrt{2}}(e_i,-e_i)$, with $i=1,2,...,n$, and let $R_{12}(\theta)$ be given by
\begin{align*}
 R_{12}(\theta) = \sum_{i=1}^n & \beta_i \ot \beta_i + \sum_{i=3}^n \gamma_i \ot \gamma_i \\
&+ \cos(\theta) \gamma_1 \ot \gamma_1 - \sin(\theta) \gamma_1 \ot \gamma_2  + \sin(\theta) \gamma_2\ot \gamma_1 + \cos(\theta) \gamma_2 \ot \gamma_2.
\end{align*}
Let $R(\theta)$ be the rotation on $\R^{nk}$ given by the matrix (\ref{matrix}) with the corresponding $R_{12}(\theta)$. From the fact that $G(R(\theta)\eta) = G(\eta)$ and $R(0) = I$ we obtain
\begin{align*}
0 = &-2\frac{\partial G(R(\theta)\eta)}{\partial \theta}|_{\theta = 0} \\
& = \left[(\eta_{12} - \eta_{22})\partial_{\eta_{11}} - (\eta_{11} - \eta_{21})\partial_{\eta_{12}} - (\eta_{12} - \eta_{22})\partial_{\eta_{21}} + (\eta_{11} - \eta_{21})\partial_{\eta_{22}}\right]G(\eta).
\end{align*}
By introducing the logarithmic derivatives $h_i' = g_i'/g_i$ the last expression becomes
\begin{equation*}
(\eta_{12} - \eta_{22})h_1'(\eta_{11}) - (\eta_{11} - \eta_{21})h_2'(\eta_{12}) - (\eta_{12} - \eta_{22}) h_1'(\eta_{21}) + (\eta_{11} - \eta_{21})h_2'(\eta_{22}) =0.
\end{equation*}
Differentiating with respect to the variable $\eta_{11}$ we obtain
\begin{equation*}
 (\eta_{12} - \eta_{22})h_1''(\eta_{11}) - h_2'(\eta_{12}) + h_2'(\eta_{22}) =0.
\end{equation*}
Finally, differentiating with respect to $\eta_{22}$ yields
\begin{equation*}
 h_1''(\eta_{11}) = h_2''(\eta_{22})\,,
\end{equation*}
and since the variables $\eta_{11}$ and $\eta_{22}$ are independent we conclude that both logarithmic second derivatives are constant. The argument above can be reproduced for $\gamma_1$ and $\gamma_j$ yielding $h_1'' = h_j'' = C$ for all $j =1,2,...,n$. This proves that all $g_i$'s are Gaussians with the same covariance, and thus $g$ will itself be a Gaussian.\\

The two last steps (reduction to the smooth non-vanishing case) plainly follows the argument of Hundertmark and Zharnitsky \cite{HZ}. This idea originally appeared in a paper by Carlen \cite{Car}. We denote by $P_{\epsilon}$ the convolution with the Gaussian kernel on $\R^{nk}$
\begin{equation*}
\varphi_{\epsilon}(\eta) = \frac{1}{(2\pi\epsilon)^{nk/2}} \,e^{-\tfrac{|\eta|^2}{2\epsilon}}\,,
\end{equation*}
and by $Q_{\epsilon}$ the convolution with the Gaussian kernel on $\R^n$
\begin{equation*}
\phi_{\epsilon}(y) = \frac{1}{(2\pi\epsilon)^{n/2}} \,e^{-\tfrac{|y|^2}{2\epsilon}}.
\end{equation*}
\\
{\bf Step 4.} Let $g \in L^p(\R^n)$ be such that $G(\eta)$ satisfies the property ($\star$). Assume $Q_{\epsilon}(g)$ never vanishes as $\epsilon \to 0$. Then $g$ is a Gaussian.\\
\\
Observe that $P_{\epsilon}(G)$ inherits the rotational symmetries of $G$, and since
\begin{equation}\label{PQ}
P_{\epsilon}(G)(\eta) = Q_{\epsilon}(g)(\eta_1)\,Q_{\epsilon}(g)(\eta_2)\,...\,Q_{\epsilon}(g)(\eta_k)\,,
\end{equation}
and $Q_{\epsilon}(g)$ is smooth and non-vanishing, we conclude by Step 3 that it must be a Gaussian. As $g \in L^p(\R^n)$, we have $g = \lim_{\epsilon \to 0} Q_{\epsilon}(g)$ and this implies that $g$, being a limit of Gaussians, is also a Gaussian.\\
\\
{\bf Step 5.} Let $g \in L^p(\R^n)$ be such that $G(\eta)$ satisfies the property ($\star$). Then $Q_{\epsilon}(g)$ never vanishes as $\epsilon \to 0$.\\
\\
Indeed, take absolute values in (\ref{PQ}) and apply the convolution operator $P_{\lambda}$ again
\begin{equation*}
P_{\lambda}|P_{\epsilon}(G)|(\eta) = Q_{\lambda}|Q_{\epsilon}(g)|(\eta_1)\,Q_{\lambda}|Q_{\epsilon}(g)|(\eta_2)\,...\,Q_{\lambda}|Q_{\epsilon}(g)|(\eta_k)\,.
\end{equation*}
Again, $P_{\lambda}|P_{\epsilon}(G)|$ inherits all the rotational symmetries of $P_{\epsilon}(G)$, in particular those of $G$. Since $Q_{\epsilon}(g) \to g$ in $L^p(\R^n)$, as $\epsilon \to 0$, we conclude that $Q_{\epsilon}(g)$ is not the zero function for small $\epsilon$. Since convolution with a Gaussian improves positivity, $Q_{\lambda}|Q_{\epsilon}(g)|$ is a strictly positive smooth function. By Step 4 we conclude that $|Q_{\epsilon}(g)|$ is a Gaussian, and thus never vanishes for small $\epsilon$.\\
\\
By putting $g = \tf$ in Steps 1-5 we are led to the conclusion that $\tf$ must be a Gaussian, and then so is $f$.


\section{Proof of Theorem \ref{thm5}: sharp cone estimates}
This final section is devoted to a brief proof of Theorem {\ref{thm5}, in which we follow the basic ideas of Foschi \cite[sections 5 and 6]{F}. Let us prove first the case $n=3$, $q=4$, which corresponds to (\ref{cone2}). From now on we shall write
\begin{equation*}
 g(|\omega|, \omega) = f(\omega)\,,
\end{equation*}
and assume that $f$ is a smooth, compactly supported function. Observe that 
\begin{equation*}
\|\widehat{g\dsigma}\|_{L^{4}(\R^{4})}^2 = \|\bigl(\widehat{g\dsigma}\bigr)^2\|_{L^{2}(\R^{4})} = \|g\dsigma * g\dsigma\|_{L^{2}(\R^{4})}\,,
\end{equation*}
where, in the case of the cone, we identify
\begin{equation*}
g\dsigma(\tau, \omega) = f(\omega) \dfrac{\delta(\tau - |\omega|)}{|\omega|}.
\end{equation*}
Therefore we can write 
\begin{equation}\label{conv g}
g\dsigma * g\dsigma (\tau, \omega) = \int_{\R^3 \times \R^3} \dfrac{f(\eta)f(\xi)}{|\eta|\,|\xi|} \,\delta_3(\omega - \eta - \xi)\, \delta(\tau - |\eta| - |\xi|) \, \deta\, \dxi\,,
\end{equation}
and we observe that $g\dsigma * g\dsigma$ is supported in the closure of the region
\begin{equation*}
C_{++} = \{(\tau, \omega) \in \R \times \R^3: \, \tau > |\omega|\}.
\end{equation*}
For each choice of $(\tau, \omega) \in C_{++}$, we denote by $\langle \cdot, \cdot \rangle_{(\tau, \omega)}$ the $L^2$-inner product associated with the measure 
\begin{equation}\label{measure delta}
\mu_{(\tau, \omega)} : = \delta_3(\omega - \eta - \xi)\, \delta(\tau - |\eta| - |\xi|) \, \deta\, \dxi\\,
\end{equation}
and by $\|\cdot\|_{(\tau,\omega)}$ the corresponding norm. From (\ref{conv g}) and Cauchy-Schwarz inequality we have
\begin{align}\label{Cauchy}
\begin{split}
g\dsigma * g\dsigma (\tau, \omega) = \langle  \frac{f(\eta)f(\xi)}{|\eta|^{1/2}\,|\xi|^{1/2}}&,\frac{1}{|\eta|^{1/2}\,|\xi|^{1/2}} \rangle_{(\tau, \omega)} \\
\\
& \leq \left\|\frac{f(\eta)f(\xi)}{|\eta|^{1/2}\,|\xi|^{1/2}}\right\|_{(\tau, \omega)} \left\|\frac{1}{|\eta|^{1/2}\,|\xi|^{1/2}}\right\|_{(\tau, \omega)}.
\end{split}
\end{align}
In \cite[Lemma 5.2]{F} it is proved that for each $(\tau, \omega) \in C_{++}$
\begin{equation}\label{lemma5.2}
\left\|\frac{1}{|\eta|^{1/2}\,|\xi|^{1/2}}\right\|_{(\tau, \omega)} = (2\pi)^{1/2}.
\end{equation}
Therefore, combining (\ref{Cauchy}) and (\ref{lemma5.2}) we obtain
\begin{align}\label{Cauchy2}
\begin{split}
\|\widehat{g\dsigma}\|_{L^{4}(\R^{4})}^4  &= \|g\dsigma * g\dsigma\|_{L^{2}(\R^{4})}^2 \leq 2\pi \int_{C_{++}}  \left\|\frac{f(\eta)f(\xi)}{|\eta|^{1/2}\,|\xi|^{1/2}}\right\|_{(\tau, \omega)}^2 \dtau \,\domega\\
& = 2\pi \int_{\R^3} \frac{|f(\eta)|^2|f(\xi)|^2}{|\eta|\,|\xi|} \deta\, \dxi = 2\pi \|g\|_{L^2(S;d\sigma)}^4\,,
\end{split}
\end{align}
and this proves (\ref{cone2}). From the Cauchy-Schwarz condition, we know that equality in (\ref{Cauchy2}) can only be attained if there is a function $F: C_{++} \to \C$ such that 
\begin{equation*}
\frac{f(\eta)f(\xi)}{|\eta|^{1/2}\,|\xi|^{1/2}} = F(\tau, \omega) \frac{1}{|\eta|^{1/2}\,|\xi|^{1/2}}\,,
\end{equation*}
for almost all $(\eta,\xi)$ (with respect to the measure (\ref{measure delta})) in the support of the measure (\ref{measure delta}), and almost all $(\tau, \omega) \in C_{++}$, with respect to the Lebesgue measure in $\R \times \R^3$. This means that 
\begin{equation}\label{charac1}
 f(\eta)f(\xi) = F(|\eta| + |\xi|, \eta + \xi)\,,
\end{equation}
for almost all $\eta, \xi \in \R^3$. The locally integrable functions $f$ satisfying property (\ref{charac1}) were characterized by Foschi in \cite[Proposition 7.23]{F} and they turn out to be
\begin{equation*}
g(|\omega|, \omega) = f(\omega) =  e^{A|\omega| + b\cdot \omega + C}\,,
\end{equation*}
where $A, C \in \C$, $b \in \C^3$ and $|\Re(b)| < -\Re(A)$ (this last condition to ensure that $g \in L^2(S;d\sigma)$).\\

The proof for the case $n=2$, $q=6$, which corresponds to (\ref{cone1}), follows exactly the same outline. Here we will have
\begin{equation*}
\|\widehat{g\dsigma}\|_{L^{6}(\R^{3})}^3 = \|\bigl(\widehat{g\dsigma}\bigr)^3\|_{L^{2}(\R^{3})} = \|g\dsigma * g\dsigma*g\dsigma\|_{L^{2}(\R^{3})}\,,
\end{equation*}
where 
\begin{equation*}
g\dsigma * g\dsigma* g\dsigma(\tau, \omega) = \int_{\R^2 \times \R^2} \dfrac{f(\eta)f(\xi)f(\zeta)}{|\eta|\,|\xi|\,|\zeta|} \,\delta_2(\omega - \eta - \xi -\zeta)\, \delta(\tau - |\eta| - |\xi|-\zeta) \, \deta \dxi\dzeta.
\end{equation*}
For each $(\tau, \omega)$ on the region $C_{+} = \{(\tau, \omega) \in \R \times \R^2: \, \tau > |\omega|\}$, consider the measure 
\begin{equation*}
\nu_{(\tau, \omega)} : = \delta_2(\omega - \eta - \xi -\zeta)\, \delta(\tau - |\eta| - |\xi|-|\zeta|) \, \deta\, \dxi\,\dzeta\\.
\end{equation*}
As in (\ref{Cauchy}), by Cauchy-Schwarz inequality
\begin{equation*}
 g\dsigma * g\dsigma* g\dsigma(\tau, \omega) \leq \left\|\frac{f(\eta)f(\xi)f(\zeta)}{|\eta|^{1/2}\,|\xi|^{1/2}\,|\zeta|^{1/2}}\right\|_{(\tau, \omega)} \left\|\frac{1}{|\eta|^{1/2}\,|\xi|^{1/2}\,|\zeta|^{1/2}}\right\|_{(\tau, \omega)}.
\end{equation*}
In \cite[Lemma 6.1]{F} it is proved that for each $(\tau, \omega) \in C_{+}$
\begin{equation}\label{lemma6.1}
\left\|\frac{1}{|\eta|^{1/2}\,|\xi|^{1/2}\,|\zeta|^{1/2}}\right\|_{(\tau, \omega)} = 2\pi.
\end{equation}
Therefore 
\begin{align}\label{Cauchy3}
\begin{split}
\|\widehat{g\dsigma}&\|_{L^{6}(\R^{3})}^6  = \|g\dsigma * g\dsigma*g\dsigma\|_{L^{2}(\R^{3})}^2 \\
\\
&\leq 4\pi^2 \int_{C_{+}}  \left\|\frac{f(\eta)f(\xi)f(\zeta)}{|\eta|^{1/2}\,|\xi|^{1/2}\,|\zeta|^{1/2}}\right\|_{(\tau, \omega)}^2 \dtau \,\domega\\
\\
& = 4\pi^2 \int_{\R^3} \frac{|f(\eta)|^2|f(\xi)|^2|f(\zeta)|^2}{|\eta|\,|\xi|\,|\zeta|} \deta\, \dxi\, \dzeta = 4\pi^2 \|g\|_{L^2(S;d\sigma)}^6\,,
\end{split}
\end{align}\\
which proves (\ref{cone1}). As in the previous case, equality happens in (\ref{Cauchy3}) if and only if there is a function $F:C_{+} \to \C$ such that
\begin{equation}\label{charac2}
f(\eta)f(\xi)f(\zeta) = F(|\eta| + |\xi|+ |\zeta|, \eta + \xi+ \zeta)\,,
\end{equation}
for almost all $\eta, \xi, \zeta \in \R^2$. The locally integrable functions $f$ satisfying (\ref{charac2}) were also characterized in \cite[Proposition 7.19]{F}, and they are
\begin{equation*}
g(|\omega|, \omega) = f(\omega) =  e^{A|\omega| + b\cdot \omega + C}\,,
\end{equation*}
where $A, C \in \C$, $b \in \C^2$ and $|\Re(b)| < -\Re(A)$. This concludes the proof.

\section*{Acknowledgments}
I am deeply grateful to my advisor William Beckner for the unceasing encouragement and all the fruitful discussions on Analysis. I am also thankful to the Clay Mathematics Institute and the National Science Foundation for supporting my participation in the Summer School on Evolution Equations 2008 in the Eidgen\"ossische Technische Hochschule - ETH, Zurich, where this work was initiated.

\end{document}